\documentclass{article}

\usepackage{arxiv}

\usepackage[utf8]{inputenc} 
\usepackage[T1]{fontenc}    
\usepackage{hyperref}       
\usepackage{url}            
\usepackage{booktabs}       
\usepackage{amsfonts}       
\usepackage{nicefrac}       
\usepackage{microtype}      
\usepackage{lipsum}		
\usepackage{graphicx}
\usepackage{natbib}
\usepackage{doi}

\title{Finite groups whose maximal subgroups have only soluble proper subgroups}

\author{ \href{https://orcid.org/0000-0003-3028-8490}{\includegraphics[scale=0.06]{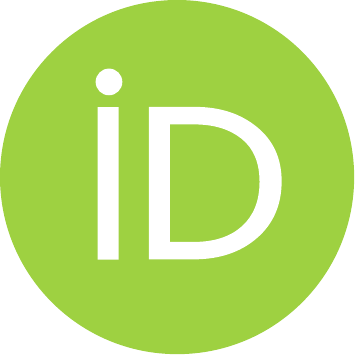}\hspace{1mm}Daria Lytkina} \\
	Sobolev Institute of Mathematics SB RAS\\
	Novosibirsk, 630090, Russia\\
	\texttt{daria.lytkin@gmail.com} \\
	{\bf Archil Zhurtov} \\
	Kabardino-Balkarian State university\\
	Nalchik, 360004, Russia\\
	\texttt{zhurtov\_a@mail.ru} \\
}

\renewcommand{\shorttitle}{Finite groups with soluble second maximal subgroups}

\hypersetup{
pdftitle={Finite groups whose maximal subgroups have only soluble proper subgroups},
pdfsubject={q-bio.NC, q-bio.QM},
pdfauthor={Daria Lytkina},
pdfkeywords={Finite group, Maximal subgroup, Solubility},
}

\begin{document}
\maketitle

\begin{abstract}
	We give a description of a finite group whose maximal subgroups possess only soluble proper subgroups, which implies the answer to the well-known question on composition factors of finite groups, whose second maximal subgroups are soluble.
\end{abstract}


\section{Introduction}
In 1968 J.~Thompson published a paper \cite{tho1968} containing classification of minimal simple groups by proving the fact that a finite simple nonabelian group, all of whose proper subgroups are soluble, is isomorphic to one of the groups of the following list:

{\bf List 1.} \emph{Minimal simple nonabelian groups}:
\begin{enumerate}
  \item $L_2(2^p)$, $p$ is a prime;
  \item $L_2(3^p)$, $p$ is an odd prime;
  \item $L_2(p)$, $p$ is a prime, $p>3$ and $p\equiv\pm 2\pmod{5}$;
  \item $Sz(2^p)$, $p$ is an odd prime;
  \item $L_3(3)$.
\end{enumerate}

NOTE. Throughout, we use the notation from \cite{atlas}.

The goal of the paper is to get a description of a finite group whose maximal subgroups possess only soluble proper subgroups.

Of course, the desired result will generalize the mentioned classification theorem by Thompson. The only difference from the original study is that we use the classification theorem on finite simple groups stating that every nonabelian finite simple group is isomorphic to one of the groups from the following

{\bf List 2.} \emph{Nonabelian finite simple groups}:
\begin{enumerate}
  \item One of 27 sporadic groups, including the Tits group $^2F_4(2)'$;
  \item $A_n$, $n\geq 5$;
  \item $L_n(q)$, $n\geq 2$, $(n,q)\notin\{(2,2),(2,3),(2,4),(2,5),(2,9),(3,2),(4,2)\}$;
  \item $U_n(q)$, $n\geq 3$, $(n,q)\neq(3,2)$;
  \item $S_n(q)$, $n\geq 4$, $n$ is even, $(n,q)\notin\{(4,2),(4,3)\}$;
  \item $O_{2n-1}(q)$, $O_{2n}^+(q)$, $O_{2n}^-(q)$, $n\geq 4$;
  \item $G_2(q)$, $q>2$;
  \item $F_4(q)$, $E_6(q)$, $E_7(q)$, $E_8(q)$;
  \item $Sz(2^{2n+1})$; $^2G_2(3^{2n+1})$; $^2F_4(2^{2n+1})$, $n\geq 1$; $^3D_4(q)$; $^2E_6(q)$.
\end{enumerate}

Note that the list is designed to avoid repetitions. For example, part 3 does not include groups $L_2(4)$ and $L_2(5)$ isomorphic with $A_5$, $L_2(9)$ isomorphic with $A_6$, $L_4(2)$ isomorphic with $A_8$ and $L_3(2)$ isomorphic with $L_2(7)$.

\section{Main result}

{\bf THEOREM.} {\it Suppose that $G$ is a finite group. Then all proper subgroups of any of its maximal subgroups are soluble if and only if one of the following conditions hold:

$(1)$ $G$ is soluble;

$(2)$ $G/\Phi(G)$ is a minimal simple group, i.\,e. it is isomorphic to one of the groups from List $1$.

$(3)$ $G$ contains a normal subgroup $G_0$ of prime index, and $G/\Phi(G_0)$ is a minimal simple group.

$(4)$ $G/\Phi(G)$ is a simple group, isomorphic to one of the groups from the following List $3$.}

{\bf List 3.} 
\begin{enumerate}
	\item $L_2(2^{rs})$, $r$ and $s$ are primes;
	\item $L_2(3^{rs})$, $r$ and $s$ are odd primes;
	\item $L_2(p)$, $p$ is a prime and $p\equiv\pm1\pmod 5$;
	\item $L_2(p^r)$, $r$ is an odd prime, $p=5$ or $p$ is a prime, such that $p\equiv\pm2\pmod5$;
	\item $L_2(9)\simeq A_6$, $U_3(3)$, $Sz(2^{rs})$, $r$ and $s$ are odd primes.
\end{enumerate}

PROOF. Suppose that case (1) is not true, i.\,e. $G$ is insoluble. Suppose first that every of its proper subgroups is soluble. Let $N$ be a maximal normal subgroup of $G$. Then $\overline{G}=G/N$ is a simple group and $N$ is soluble. Therefore, $\overline{G}$ is nonabelian, all of its proper subgroups are soluble, and $\overline{G}$ is isomorphic to one of the groups from List 1. If $M$ is a maximal subgroup of $G$, then $MN\neq G$ and hence $N\le M$. So, $N$ is contained in the intersection $\Phi(G)$ of maximal subgroups and case (2) is true.

Next suppose that $G$ has an insoluble maximal subgroup $M$. If $M\trianglelefteq G$, then $G/M$ is a group of prime order, all proper subgroups of $M$ are soluble, and case (3) is true.

And if $G$ has no proper normal insoluble subgroup and $N$ is a maximal in $G$ normal soluble subgroup, then $E=G/N$ is a simple nonabelian group and among its maximal subgroups there is an insoluble subgroup $M$. By condition all proper subgroups in $M$ are soluble and, therefore, $M/\Phi(M)$ is a simple group isomorphic to one of the groups from List 1.

If $MN=G$, then $G/N\simeq M/N\cap M$. Since all proper subgroups of $M$ are soluble and $G$ has exactly one insoluble composition factor, then all proper subgroups of $G$ are soluble and, using what we proved before, $G/\Phi(G)$ is a simple group from List 1. So, $N$ lies in any maximal subgroup of the group $G$, i.e. $N\le\Phi(G)$.

Now, to prove the Theorem, it is suffice to confirm the following statement.

{\bf Proposition.} {\it Let $G$ be a finite simple group containing a proper insoluble subgroup, and let every such insoluble subgroup modulo its Frattini subgroup is isomorphic to a group from List $1$. Then $G$ is isomorphic to one of the groups from List $3$.}

Since $G$ is isomorphic to one of the groups from List 2, we have to inspect this list.

Information about maximal subgroups of sporadic groups provided in \cite{atlas} shows that none of these groups fits to be the group $G$.

Suppose that $G=L_n(q)$, $q=p^s$, $p$ is a prime, $s\in\mathbb{N}$. If $n=2$, then by \cite[Satz II,8,27]{hup} $G$ is isomorphic to one of the following groups:
\begin{itemize}
	\item $L_2(2^{rs})$, $r$, $s$ are primes;
	\item $L_2(3^{rs})$, $r$, $s$ are odd primes;
	\item $L_2(p)$, $p\equiv\pm1 \pmod 5$;
	\item $L_2(p^r)$, $p=5$ or $p\equiv\pm2 \pmod 5$, $r$ is an odd prime;
	\item $L_2(9)\simeq A_6$.
\end{itemize}
Note, along the way, that none of the alternating groups $A_m$ for $m\ge 7$ satisfies the conditions of the Theorem.

If $n\ge3$, then the stabilizer $S$ in $G$ of a one-dimensional subspace of projective space $P_{n-1}$ under natural action of $G$ is isomorphic to $\left([q^{n-1}]:GL_{n-1}\right)/(q-1,n)$. If $S$ is insoluble, then it does not satisfy the conditions of the Proposition. It is soluble given $n=2$ and $q=2$ or $3$. For these two cases $G$ is in List 1.

Another case is when $G=S_n(q)$, $n\ge 4$, $n$ is even, $(n,q)\notin\{(4,2),(4,3)\}$. Then the maximal subgroup of the smallest index in $G$ is insoluble by \cite{maz1993e} and does not satisfy the conditions of the Proposition.

Let $G=U_n(q)$. 

If $n=3$, then Tables 8.5 and 8.6 from \cite{brh} show that the conditions of the Proposition are met only by $G=U_3(3)$.

If $n=4$, then one of the maximal 2-local subgroups of the group $U_4(q)$, which is isomorphic to $\left([q^4]:SL_2(q^2):(q-1)\right)/(q+1,4)$ is insoluble and does not satisfy the conditions of the Proposition.

Now, if $n\ge 5$, then the monomial subgroup of $U_n(q)$ isomorphic to $\left((q+1)^{n-1}/(q+1,n)\right):S_n$ is insoluble and does not meet the conditions of the Proposition.

If $G$ is an orthogonal group in dimension at least 7, its subgroup of the smallest index is insoluble and does not satisfy the conditions of the Proposition (see \cite{vama1994e}).

Same is true for exceptional groups of types $F_4$, $E_6$, $E_7$, and $E_8$, and also for groups of the type $^2E_6$ (see \cite{va1996e,va1997e,va1998e}). Maximal subgroups of simple groups of the type $G_2$, $Sz(2^{2n+1})$, $^2G_2(3^{2n+1})$, $^2F_4(2^{2^{2n+1}})$, and $^3D_4(q)$ are listed in \cite{brh} (Tables 8.16, 8.30, 8.41, 8.42, 8.43, 8.51). Except for the groups $Sz(2^{2n+1})$, none of the groups in those tables satisfies the conditions of the Proposition.

For Suzuki groups, the conditions of the Proposition are met only by the groups $Sz(2^{rs})$, where $r$ and $s$ are odd primes. This proves the Proposition along with the Theorem.

The Theorem ovbiously implies the answer to the well-known question on composition factors of finite groups, whose second maximal subgroups are soluble.

{\bf Corollary.} {\it Consider a finite group $G$, such that every maximal subgroup of its maximal subgroup is soluble. If $C$ is a composition factor of $G$, then $C$ is either a group of prime order, or $C$ is isomorphic to one of the groups from List $1$ or List $3$.}

\end{document}